\input amstex
\input xy
\xyoption{all}
\loadeufm
\loadeusm
\loadeusb
\loadeurb
\loadeurm

\documentstyle{amsppt}
\mag1200\vsize.833\vsize

\NoRunningHeads
\NoBlackBoxes
\topmatter

\title
PURITY FOR SIMILARITY FACTORS
\endtitle
\author
Ivan Panin
\endauthor

\date  June 2003
\enddate

\abstract Let $R$ be a regular local ring, $K$ its field of
fractions and $A_1$, $A_2$ two Azumaya algebras with involutions
over $R$. We show that if $A_1 \otimes_R K$ and $A_1 \otimes_R K$
are isomorphic over $K$, then $A_1$ and $A_2$ are isomorphic over
$R$. In particular,
if two quadratic spaces over the ring $R$ become
similar over $K$ then these two spaces are
similar already over $R$. The results are consequences of
a purity theorem for similarity factors.
\endabstract


\thanks
The author thanks very much for the support the RTN-Network
HPRN-CT-2002-00287 , the grant of the years 2001 - 2003 of
the Russian Science Support Foundation at the Russian Academy of
Science, the grant INTAS-99-00817 and the University of Lausanne.
\endthanks

\endtopmatter

\document

\head Introduction
\endhead
Let $R$ be a regular local ring, $K$ its field of fractions. Let
$(A_1, \sigma_1)$ and $(A_2, \sigma_2)$ be two Azumaya algebras
with involutions over $R$ (see right below for a precise
definition). Assume that
$(A_1, \sigma_1)\otimes_R K$ and
$(A_2,\sigma_2)\otimes_R K$
are isomorphic. Are
$(A_1, \sigma_1)$ and
$(A_2, \sigma_2)$ isomorphic too?
We show that this is true if $R$
is a regular local ring containing a field of characteristic
different from $2$.
If $A_1$ and $A_2$ are both the $n \times n$ matrix algebra
over $R$ and the involutions are symmetric
then $\sigma_1$ and $\sigma_2$ define
two quadratic spaces $q_1$ and $q_2$ over $R$ up to similarity factors.
In this particular case the result looks as follows:
if $q_1 \otimes_R K$ and $q_1 \otimes_R K$ are similar then
$q_1$ and $q_2$ are similar too.

Grothendieck [G] conjectured that, for any reductive group scheme
$G$ over $R$, rationally trivial $G$-homogeneous spaces are
trivial. Our result corresponds to the case  when  $G$ is the
projective unitary group ${\roman P \roman U}_{A,\sigma}$ for an
Azumaya algebra with involution over $R$. If $R$ is an essentially
smooth local $k$-algebra and $G$ is defined over $k$ (we say that
$G$ {\it is constant}) Grothendieck's conjecture has been proved
in most cases: by Colliot-Th\'el\`ene and Ojanguren \cite{\bf
C-TO} for a perfect infinite field $k$ and then by Raghunathan
\cite{\bf R} for any infinite $k$. One notable open case is that
of a finite base field. For a non-constant group $G$ only few
cases have been proved: when $G$ is a torus, by
Colliot-Th\'el\`ene and Sansuc \cite{\bf C-TS}, when $G$ is the
group $\roman{SL}_1(D)$ of norm one elements of an Azumaya
$R$-algebra $D$, by Panin and Suslin \cite{\bf PS}, when $G$ is
the unitary group $\roman U_{A, \sigma}$, by Panin and Ojanguren
\cite{\bf Oj-P1}, when $G$ is the special unitary group $\roman S
\roman U_{A, \sigma}$, by Zainoulline \cite{\bf Z}.
Recall as well that for semi-simple group schemes $G$ over a discrete
valuation ring the conjecture has been proved by Nisnevich in
\cite{\bf Ni}.

The paper is organized as follows.
Section 1 contains a reduction of the main theorem (Th. 1.1)
to a purity theorem for similarity factors (Th. 1.3).
Section 2 is devoted to a theorem of Nisnevich and its Corollaries.
The rest of the text is devoted to the proof of Theorem 1.3. The proof is
given in \S8. It is based on the Specialization Lemma (stated in \S3 and proved
in \S4), the Equating Lemma (\S5) and the Unramifiedness Lemma (\S5).

\smallskip
The author thanks very much M.Ojanguren for a lot of useful discussions
on the subject of the present article.

\head\S1.
Rationally isomorphic Azumaya algebras with involutions
are locally isomorphic
\endhead

Let $R$
be a regular local ring containing a field $k$
$({\text{char}}(k)\ne2)$
and let
$K$
be its quotient field. By an $R$- Azumaya algebra
with involution $(A, \sigma)$  we mean (see \cite{\bf Oj-P1})
an $R$-algebra $A$ which is an Azumaya algebra over its center $Z(A)$ equipped with  an involution $\sigma: A \to A^{op}$,
such that  $Z(A)$ is either $R$ itself or
an \'etale quadratic extension of $R$ such that $Z(A)^\sigma=R$.
\proclaim{1.1. Theorem (Main)}
Let
$(A_1, \sigma_1)$
and
$(A_2, \sigma_2)$
be two Azumaya algebras with involutions over the ring $R$.
If the Azumaya algebras with involutions
$(A_{1,K}, \sigma_{1,K})$
and
$(A_{2,K}, \sigma_{2,K})$
are isomorphic,
then
$(A_1, \sigma_1)$
and
$(A_2, \sigma_2)$
are already isomorphic.
\endproclaim
\subhead
Reduction to a Purity Theorem
\endsubhead
Since
$
(A_{1,K}, \sigma_{1,K})\simeq (A_{2,K}, \sigma_{2,K})
$,
one concludes that the two Azumaya algebras
$
A_{1,K}
$
and
$
A_{2,K}
$
over
$
Z_K
$
are isomorphic. Thus
$
A_1\simeq A_2
$
($
Z
$
is regular semilocal as an \'etale quadratic extension of
$
R
$).
Therefore one may assume that
$
A_1=A_2
$,
$
Z_1=Z_2
$
and we have two involutions
$
\sigma_1
$
and
$
\sigma_2
$
on the same algebra
($
A
$
over
$
Z
$
and
$
Z
$
is a quadratic \'etale extension of
$
R
$
(or
$
Z
$
just coincides with
$
R
$).

Now consider the composite
$
A@>\sigma_2>>A^{op}@>\sigma_1^{-1}>>A
$.
It is an Azumaya algebra isomorphism. Thus it is of the form
$
\text{Int}(\alpha)
$
for an element
$
\alpha \in A^*
$.
Thus
$
\sigma_1\circ \text{Int}(\alpha)=\sigma_2
$
and
$
\alpha
$
is symmetric with respect to
$
\sigma_1
$.
Therefore we have two hermitian spaces over
$
(A, \sigma_1)
$,
namely
$
(A, 1)
$
and
$
(A, \alpha)
$.
Set
$$
h_1=(A, 1)\text{ and }h_2=(A, \alpha).
$$
Since
$
(A_K, \sigma_{1,K})
$
is isomorphic to
$
(A_K, \sigma_{2,K})
$,
$
h_{1,K}
$
is similar to
$
h_{2,K}
$
i.e. there exist an element
$
a\in K^*
$
and an isometry
$
a \cdot h_{1,K}\simeq h_{2,K}
$.

We will prove (this suffices to prove the theorem) that
there exists an element
$
b\in R^*
$
such that
$
b \cdot h_1\simeq h_2
$
over
$
(A, \sigma_1)
$.
To find the desired element
$
b\in R^*
$,
it suffices to find a similarity
factor
$
b_\eufm m \in K^*
$
of the space
$
h_{1, K}
$
and a unit
$
a_\eufm m \in R^*
$
such that
$
a=a_\eufm m \cdot b_\eufm m
$.
In fact, if
$
b_\eufm m \in K^*
$,
$
a_\eufm m \in R^*
$
are the mentioned elements then one has a chain of relations
$
(h_2\perp -a_{\eufm m} \cdot h_1)_K\simeq h_{2,K}\perp -a_\eufm m \cdot b_\eufm m \cdot h_{1,K}=h_{2,K}\perp -a \cdot h_{1,K}\simeq h_{2,K}\perp -h_{2,K}
$.
Thus
$
(h_2\perp -a_{\eufm m} \cdot h_1)_K
$
is hyperbolic and
by the main theorem of \cite{\bf{Oj-P1}} the space
$
h_2\perp -a_{\eufm m} \cdot h_1
$
is hyperbolic,
whence
$
h_2\simeq a_{\eufm m}\cdot h_1
$.
Therefore putting
$
b=a_{\eufm m}
$
we get
$
h_2\simeq b \cdot h_1
$
over
$
(A, \sigma_1)
$.
It remains to find a similarity factor
$
b_{\eufm m}
$
of
$
h_{1,K}
$
and a unit
$
a_{\eufm m}\in R^*
$
such that
$
a=a_{\eufm m}\cdot b_{\eufm m}
$.
By the corollary of a theorem of Nisnevich below (\S3, Cor. 3.2), for a height
one prime ideal
$
\eufm p
$
in
$
R
$
there exist elements
$
b_{\eufm p}\in K^*
$
and
$
a_{\eufm p}\in R^*
$
such that
\roster
\item
$
b_{\eufm p}\text{ is a similarity factor of the space }h_{1,K}\text{ and}
$
\item
$
a=a_{\eufm p}\cdot b_{\eufm p}
$.
\endroster
Thus by the Purity Theorem (Theorem 1.2)  there exist a similarity factor
$
b_{\eufm m}
$
of
$
h_{1,K}
$
and a unit
$
a_{\eufm m}\in R^*
$
with
$
a=a_{\eufm m}\cdot b_{\eufm m}
$.
So we have reduced Theorem 1.1 to the Purity Theorem.
\qed
\proclaim{1.2. Theorem (Purity Theorem)}
Let
$
R
$,
$
K
$
be  as in Theorem 1.1. Let
$
(A,\sigma)
$
be an Azumaya algebra with involution over
$
R
$
and let
$
h
$
be the hermitian space
$
(A,1)
$
over
$
(A,\sigma)
$.
Let
$
a\in K^*
$.
Suppose that for each prime ideal of height $1$
$
\eufm p
$
in
$
R
$
there exist
$
a_{\eufm p}\in R_{\eufm p}^*
$,
$
b_{\eufm p}\in K^*
$
with
$
a=a_{\eufm p}\cdot b_{\eufm p}
$
and
$
h_K\simeq b_{\eufm p}\cdot h_K
$.
Then there exist
$
b_{\eufm m}\in K^*
$,
$
a_{\eufm m}\in R^*
$
such that
\roster
\item
$
b_{\eufm m}\text{ is a similarity factor of the space }h_K\text{,}
$
\item
$
a=a_{\eufm m}\cdot b_{\eufm m}
$.
\endroster
\endproclaim
It is convenient for the proof to restate Theorem 1.2  it in a slightly
more technical form. For that consider the similitude group scheme
$G=\text{Sim}_{A,\sigma}$ of the Azumaya algebra with involution
$(A,\sigma)$. Recall that for an $R$-algebra $S$ the $S$-points of
$G$ are those
$\alpha \in (A \otimes_R S)^*$
for which
$\alpha^\sigma\cdot\alpha \in S^*$.
Further consider a group scheme morphism
$\mu: G \to \Bbb G_m$
which takes a similitude
$\alpha \in G(S)$
to its similarity factor
$\mu (\alpha)=\alpha^\sigma\cdot\alpha \in S^*$.
Finally for an $R$-algebra $S$ consider the group
$\Cal F(S)=S^*/ \mu (G(S))$.
For an element $a \in S^*$ we will often write
$\bar a$ for its class in $\Cal F(S)$.

\proclaim{1.3. Theorem }
Let $R$, $K$ and $(A,\sigma)$ be  as in Theorem 1.2.
Let $a \in K^*$. If for each height $1$ prime $\eufm p$ in $R$ the class
$\bar a \in \Cal F(K)$
can be lifted in
$\Cal F(R_{\eufm p})$,
then $\bar a$ can be lifted in $\Cal F(R)$.
\endproclaim

\remark
{Remark}
Theorems 1.2 and 1.3 are equivalent. In fact,
the group
$\mu (G(R))$
coincides with the group $G_R(h)$
of similarity factors of the hermitian space
$h=(A,1)$.
\endremark

\remark {Remark} It is quite plausible that the method of \cite{\bf
Z1} could be adapted to prove Theorem 1.3.
\endremark

\head\S2.
A theorem of Nisnevich
\endhead
Let
$
R
$
be a discrete valuation ring containing a field and let
$
K
$
be its quotient field.
Let
$
(A, \sigma)
$
be an Azumaya algebra with involution over
$
R
$.
The following theorem is a consequence of a theorem of Nisnevich on principal $G$-bundles.
(\cite{\bf Ni}, Theorem ??).
\proclaim{2.1. Theorem (Nisnevich)}
Let
$
h_1
$
and
$
h_2
$
be two hermitian spaces over
$
(A, \sigma)
$.
Suppose
$
h_{1,K}
$
is similar to
$
h_{2,K}
$,
then
$
h_1
$
is similar to
$
h_2
$.
\endproclaim
This Theorem is a particular case of the theorem of Nisnevich just mentioned,
namely the case when $G$ is  the projective unitary group scheme
$PU_{h_1}$ over $R$.

\proclaim{2.2. Corollary}
Let
$
h_1
$,
$
h_2
$
be two hermitian spaces over
$
(A,\sigma)
$.
Let
$
a\in K^*
$
be such that
$
h_{2,K}\simeq a\cdot h_{1,K}
$.
Then there exist an element
$
b^{\prime}\in K^*
$
and a unit
$
a^{\prime}\in R^*
$
such that
\roster
\item
$
b^{\prime}
$
is a similarity factor of the space
$
h_{1,K}
$,
\item
$
a=a^{\prime}\cdot b^{\prime}
$.
\endroster
\endproclaim
\demo{Proof}
By the theorem there exists a unit
$
a^{\prime}\in R^*
$
such that
$
a^{\prime}\cdot h_2\simeq h_1
$.
Thus one has a chain of relations
$$
a\cdot (a^{\prime})^{-1}\cdot h_{1,K}\simeq a\cdot h_{2,K}\simeq a^2 \cdot h_{1,K}\simeq h_{1,K}\ .
$$
Therefore
$
b^{\prime}=a\cdot (a^{\prime})^{-1}
$
is a similarity factor of the space
$
h_{1,K}
$
and
$
a=a^{\prime}\cdot b^{\prime}
$.
\qed
\enddemo

\proclaim{2.3. Corollary}
The kernel of the map
$H^1(R, \text{Sim}_{A,\sigma}) \to H^1(K, \text{Sim}_{A,\sigma})$
is trivial.
\endproclaim

\demo{Proof}
The group scheme
$\text{Sim}_{A,\sigma}$
fits in an exact sequence of algebraic groups
$$
0 \to R_{Z / R}(\Bbb G_{m, Z}) \to \text{Sim}_{A,\sigma} \to \text{PU}_{A,\sigma} \to 0
$$
where
$R_{Z / R}(\Bbb G_{m, Z})$
is the Weil restriction of the multiplicative group
$\Bbb G_{m, Z}$.
By Hilbert's Theorem $90$
$H^1(R, R_{Z / R}(\Bbb G_{m, Z})) = H^1(Z, \Bbb G_{m, Z})) =0$.
Thus the kernel of the map
$H^1(R, \text{Sim}_{A,\sigma}) \to H^1(R, \text{PU}_{A,\sigma})$
is trivial. On the other hand the kernel of the map
$H^1(R, \text{PU}_{A,\sigma}) \to H^1(K, \text{PU}_{A,\sigma})$
is trivial by Theorem 3.1. Thus the kernel of the map
$H^1(R, \text{Sim}_{A,\sigma}) \to H^1(K, \text{Sim}_{A,\sigma})$
is trivial as well,  whence the Corollary.
\enddemo

\head\S3.
A Specialization Lemma
\endhead
In this section we state a theorem which is one of the main
ingredient in the proof of purity. The theorem itself will
be proved in \S 5 below.

Let $ k $ be a field $ ({\text{char}}(k)\neq 2) $ and let $
(A,\sigma) $ be an Azumaya algebra with involution over $ k $ (see
Section $1$ for the definition). Let $ G={\text{Sim}}_{A,\sigma} $
be the similitude group of $ (A,\sigma) $ (see the end of Section
1 for the definition), and let $
\mu:{\text{Sim}}_{A,\sigma}@>>>\Bbb G_m $ be a group morphism
which takes a similitude $ \alpha $ to its similarity factor $
\mu(\alpha)={\alpha}^{\sigma}\cdot \alpha $. The group $ G $
coincides with the similitude group of the hermitian space $
(A,1)=h $.
\definition
{3.1. Notation}
For a commutative
$
k
$-algebra
$
S
$,
set
$
{\Cal F}(S)=S^*/{\mu}(G(S))
$.
For an element
$
u\in S^*
$
we shall write in this section
$
\overline u
$
for the image of
$
u
$
in
$
{\Cal F}(S)
$.
Observe that
$
{\mu}(G(S))=G_S(h\otimes_k S)
$
is the group of similarity factors of the hermitian space
$
h\otimes_k S
$.
\enddefinition

Let
$
S
$
be a
$
k
$-algebra which is a Dedekind domain and let
$
L
$
be the quotient field of
$
S
$.
Let
$
\eufm p\subseteq S
$
be a non-zero prime ideal in
$
S
$
and let $
S_\eufm p
$
be the corresponding local ring.
\definition{3.2. Definition}
Let
$
a\in L^*
$.
The element
$
\overline a \in {\Cal F}(L)
$
is said to be {\it unramified} at a prime
$
\eufm p
$
if
$
\overline a
$
belongs to the image of the group
$
{\Cal F}(S_\eufm p)
$
in
$
{\Cal F}(L)
$.
In other terms, the element
$
\overline a
$
is  {\it unramified} at
$
\eufm p
$
if
$
a=a_{\eufm p}\cdot b_{\eufm p}
$
for certain elements
$
a_\eufm p \in S_{\eufm p}^*
$
and
$
b_{\eufm p}\in {\mu(G(L))}
$.
We denote by
$
{\Cal F}_{un}(S)
$
the subgroup in
$
{\Cal F}(L)
$
consisting of all those elements in
$
{\Cal F}(L)
$
which are unramified at each non-zero prime
$
\eufm p
$
in $
S
$.
Elements of $\Cal F(S)$ are called
$S$-{\it unramified} elements.
\enddefinition

Let $ S\supseteq k[s] $ be a finite extension of the polynomial
ring in one variable. Suppose $ S $ is a Dedekind domain. Let $
S_1=S/(s-1)S $ and $S_0=S/J $, and let $ \epsilon :S@>>>k $ be an
augmentation such that $ S/tS=S/{\text{Ker}}(\epsilon)\times
S/J=k\times S/J $ for certain ideal $ J $ in $ S $. For an element
$ v\in S $ we will write $ v_1 $ and $ v_0 $ for its images in $
S_1 $ and $ S_0 $ respectively. If furthermore $ g\in S $ be an
element coprime as to $ (s-1) $ so to $ (s) $, then the canonical
map $ S@>>>S_i $ is factorized as the composite $ S@>>>S_g@>>>S_i
$. In this case for an element $ v\in S_g $ we will write $ v_1 $
and $ v_0 $ for its images in $ S_1 $ and $ S_0 $ respectively. We
will denote $ N_{S_i/k}:S_i^*@>>>k^* $ the norm map.

\proclaim{3.3. Theorem (Specialization Lemma)} Let $ S\supseteq
k[s] $ be an integral extension of the polynomial ring in one
variable $ k[s] $ and suppose $ S $ is a Dedekind domain and $
L $  its quotient field. Let $ \eurm f \in S $ be an element
coprime  to  $ s $ and  $ (s-1) $. Let
$ u\in S_{\eurm f}^* $ be a unit. Suppose the element $ \overline
u \in {\Cal F}(L) $ is $S$-unramified, i.e. $ \overline u $
belongs to the subgroup $ \Cal F_{un}(S) $. Then the following
relation holds in the group $ \Cal F(k) $
$$
{\overline {{\epsilon}(u)}}={\overline {N_{S_1/k}(u_1)}}\cdot
{\overline {N_{S_0/k}(u_0)^{-1}}} \quad.
\tag{*}
$$
\endproclaim
\remark
{3.4. Remark}
This theorem is proved in \S 4 below. Now observe only that if
$
u\in S^*
$,
then
$
N_{S/k[s]}(u)\in k[s]^*=k^*
$
and
already
$\epsilon(u)=N_{S_1/k}(u_1)\cdot N_{S_0/k}(u_0)^{-1}$.
So there is nothing to prove in this case.
The trouble is that we do not assume
$
u\in S^*
$.
\endremark

\head\S4.
Proof of Specialization Lemma
\endhead
Let $ k $ be a field  of characteristic different of 2 and let $
(A,\sigma) $ be an Azumaya algebra with involution over $ k $ (see
Section 1 for the definition). Let $ h $ be the hermitian space $
(A,1) $. We preserve in this section notation of \S3.

Let
$
K
$
be a function field of an irreducible curve over
$
k
$
and let
$
L\supseteq K
$
be a finite field extension (separable). We will consider
in this section discrete valuations of
$
K
$
and
$
L
$
which are trivial on
$
k
$
and they will be called valuations. For valuations
$
x:K^* @>>> \Bbb Z
$
and
$
y: L^* @>>> \Bbb Z
$,
we write
$
y/x
$
if
$
y
$
extends
$
x
$.
We will need  completions to avoid dealing with
semi-local Dedekind domains.

\definition {4.1. Notation}
Let
$
y
$
be a valuation of
$
L
$.
Denote by
$
\hat L_y
$
the completion of
$
L
$
with respect to
$
y
$.
Denote by
$
\Cal O_y
$
the ring of integers associated with
$
y
$, i.e.
$
\Cal O_y=\{a\in L\ \vert\ y(a)\geqslant 0\}
$.
And denote by
$
\hat{\Cal O}_y
$
the ring of
$
y
$-integers in
$
\hat {L}_y
$, i.e.
$
\hat {\Cal O}_y=\{a\in \hat L_y\ \vert\ y(a)\geqslant 0\}
$.
We shall write
$
k(y)
$
for
the residue field of
$
y
$, i.e.
$
k(y)=\Cal O_y/ \eufm m_y=\hat {\Cal O}_y/ \hat {\eufm m}_y
$.
\enddefinition

If
$
x
$
and
$
y
$
are valuations of
$
K
$
and
$
L
$
respectively
and
$
y
$
extends
$
x
$,
then
$
\Cal O_y \supseteq \Cal O_x
$
and
$
\hat {\Cal O}_y \supseteq \hat {\Cal O}_x
$
and the ring extension
$
\hat {\Cal O}_y \supseteq \hat {\Cal O}_x
$
is integral. Thus one has norm mappings
$
N_{{\Cal O}_y / \hat {\Cal O}_x}:\hat {\Cal O}_y^* @>>> \hat {\Cal O}_x^*
$
and
$
N_{L_y / \hat K_x}:\hat L_y^* @>>> \hat K_x^*
$
(we will use below a short notation
$
N_{y/x}
$
for both of these maps). There is the norm map
$
N_{k(y)/k(x)}:k(y)^* @>>> k(x)^*
$
and two diagrams commute
$$
\CD
\hat{\Cal O}_y^* @>>> \hat {L}_y^*  \\
@VN_{y/x}VV          @VVN_{y/x}V    \\
\hat{\Cal O}_x^* @>>> \hat {K}_x^* @. , \\
\endCD
\qquad
\qquad
\qquad
\CD
\hat{\Cal O}_y^* @>>> k(y)^*  \\
@VN_{y/x}VV           @VVN_{k(y)/k(x)}^{i(y/x)}V \\
\hat{\Cal O}_x^* @>>> k(x)^* @. , \\
\endCD
$$
where
$
i(y/x)=\text{ the length of }\Cal O_y/\eufm m_x\Cal O_y
$
is the ramification index of
$
y
$
over
$
x
$.
\remark
{4.2. Remark}
Let
$
U_{A,\sigma}
$
be the unitary group of the form
$
h
$.
It is an algebraic group over
$
k
$
such that for any
$
k
$-algebra
$
R
$
the group of its
$
R
$-points
is the group
$
\{ \alpha \in (A\otimes_k R)^*\ \vert\ {\alpha}^{\sigma} \cdot \alpha = 1 \}
$.
With the notation of \S3 the group
$U_{A,\sigma}$
fits in an exact sequence of algebraic groups
$
1@>>>U_{A,\sigma}@>>>{\text{Sim}}_{A,\sigma}@>\mu>>{\Bbb G}_m@>>>1.
$
This sequence
of algebraic groups induces exact sequences of pointed sets
($R
$
is a domain,
$
K
$
is its quotient field)
$$
\CD
1@>>>{\Cal F}(R)@>\partial>>H_{et}^1(R,U_{A,\sigma})@>>>H_{et}^1(R,{\text{Sim}}_{A,\sigma}) \\
@.   @VV\theta_1V    @VV\theta_2V             @VV\theta_3V            \\
1@>>>{\Cal F}({K})@>\partial>>H_{et}^1({K},U_{A,\sigma})@>>>H_{et}^1({K},{\text{Sim}}_{A,\sigma}) . \\
\endCD
$$
In the case of a Dedekind local ring
$R$ and its quotient field $K$ the maps
$\theta_2$
and
$\theta_3$
have trivial kernels as well. This holds for
$\theta_2$
by Corollary 2.3 and for
$\theta_3$
by
\cite{\bf{Oj}}.
In particular,  in this case the map
$
\theta_1:\Cal F(R) @>>> \Cal F({K})
$
has the trivial kernel and thus it
is injective. Observe as well that
for a field
$
K
$
the map
$
{\Cal F}({K})@>\partial>>H_{et}^1({K},U_{A,\sigma})
$
is injective, i.e.
$
(\partial(a)=\partial(b)) \implies a=b
$.
\endremark

\definition
{4.3. Notation}
Let
$
y
$
be a valuation of
$
L
$
and let
$
i:L @>>> \hat L_y
$
be the inclusion. Then by Remark 4.3 the map
$
\Cal F(\hat{\Cal O}_y)@>i_*>>\Cal F(\hat{L}_y)
$
is injective and we will identify
$
\Cal F(\hat{\Cal O}_y)
$
with its image under this map. Set
$$
\Cal F_y(L)=i_*^{-1}(\Cal F(\hat{\Cal O}_y)).
$$
\enddefinition

The inclusions
$
\Cal O_y\hookrightarrow L
$
and
$
\Cal O_y\hookrightarrow \hat{\Cal O}_y
$
induce a map
$
\Cal F({\Cal O}_y) @>>> \Cal F_y(L)
$
which is injective by Remark 4.3. Both groups are subgroups of
$
\Cal F(L)
$.
The following lemma shows that
$
\Cal F_y(L)
$
coincides with the subgroup of $
\Cal F(L)
$ consisting of all elements {\it unramified} at
$
y
$.

\proclaim{4.4. Lemma}
$
\Cal F({\Cal O}_y)=\Cal F_y(L).
$
\endproclaim
\demo{Proof}
We only have to check the inclusion
$
\Cal F_y(L) \subseteq \Cal F({\Cal O}_y)
$.
Let
$
a_y \in \Cal F_y(L)
$
be an element. It determines the elements
$
a \in \Cal F(L)
$
and
$
\hat a \in \Cal F(\hat{\Cal O}_y)
$
which coincide when regarded as elements of
$
\Cal F(\hat{L}_y)
$.
We denote this common element in
$
\Cal F(\hat{L}_y)
$
 by
$
\hat a_y
$.
Let
$
\xi=\partial(a)\in H_{et}^1(L,U_{A,\sigma})
$,
$
\hat \xi=\partial(\hat a)\in H_{et}^1(\hat{\Cal O}_y,U_{A,\sigma})
$
and
$
\hat \xi_y=\partial(\hat a_y)\in H_{et}^1(\hat L_y,U_{A,\sigma})
$.
Clearly,
$
\hat \xi
$
and
$
\xi
$
both coincide with
$
\hat \xi_y
$
when regarded as elements of
$
H_{et}^1({\hat L}_y,U_{A,\sigma})
$.
Thus one can glue
$
\xi
$
and
$
\hat \xi
$
to get a
$
\xi_y \in H_{et}^1({\Cal O}_y,U_{A,\sigma})
$
which maps to $
\xi
$
under the map induced by the inclusion
$
\Cal O_y \hookrightarrow L
$
and maps to $
\hat \xi
$
under the map induced by the inclusion
$
\Cal O_y \hookrightarrow \hat{\Cal O}_y
$.

We now show that
$
\xi _y
$
has the form
$
\partial (a_y^\prime)
$
for a certain
$
a_y^\prime \in \Cal F(\Cal O_y)
$.
In fact, observe that the image
$
\zeta
$
of
$
\xi
$
in
$
H_{et}^1(L, {\text{Sim}}_{A,\sigma})
$
is  trivia. As mentioned in Remark 4.3 the map
$
H_{et}^1({\Cal O}_y,{\text{Sim}}_{A,\sigma}) @>>> H_{et}^1(L,{\text{Sim}}_{A,\sigma})
$
has the trivial kernel. Therefore the image
$
\zeta_y
$
of
$
\xi_y
$
in
$
H_{et}^1(\Cal O_y, {\text{Sim}}_{A,\sigma})
$
is trivial  as well. Thus there exists an element
$
a_y^\prime \in \Cal F(\Cal O_y)
$
with
$
\partial(a_y^\prime)=\xi_y \in H_{et}^1(\Cal O_y,U).
$

We now prove that
$
a_y^\prime
$
coincides with
$
a_y
$
in
$
\Cal F_y(L)
$.
Since
$
\Cal F(\Cal O_y)
$
and
$
\Cal F_y(L)
$
are both subgroups of
$
\Cal F(L)
$,
it suffices to show that
$
a_y^\prime
$
coincides with the element
$
a
$
in
$
\Cal F(L)
$.
By Remark 4.3 the map
$
\Cal F(L)@>\partial>> H_{et}^1(L,U_{A,\sigma})
$
is \underbar {injective}. Thus it suffices to check that
$
\partial(a_y^\prime)=\partial(a)
$
in
$
H_{et}^1(L,U_{A,\sigma})
$.
This is indeed the case because
$
\partial(a_y^\prime)=\xi_y
$
and
$
\partial(a)=\xi
$,
and
$
\xi_y
$
coincides with
$
\xi
$
when regarded over
$
L
$.
We have proved that
$
a_y^\prime \in \Cal F(\Cal O_y)
$
coincides with
$
a_y
$
in
$
\Cal F_y(L)
$.
Thus the inclusion
$
\Cal F_y(L) \subseteq \Cal F({\Cal O}_y)
$
is proved, whence the lemma.
\qed
\enddemo
\definition{4.5. Definition}
Let
$
y
$
be a valuation of
$
L
$.
Define a specialization map
$$
s(y):\Cal F_y(L)@>>>\Cal F(k(y))
$$
as the composite
$
\Cal F_y(L) @>>> \Cal F(\hat {\Cal O}_y) @>\text{res}_y>> \Cal F(k(y))
$
of the map
$
\Cal F_y(L) @>>> \Cal F(\hat {\Cal O}_y)
$
induced by the map
$
i_*
$
(see 4.4) and the map
$
\Cal F(\hat {\Cal O}_y) @>>> \Cal F(k(y))
$
induced by the residue map
$
\hat {\Cal O}_y @>>> k(y)
$.
(If we identify
$
\Cal F_y(L)
$
with
$
\Cal F({\Cal O}_y)
$
by Lemma 4.5, then the map
$
s(y):\Cal F(\Cal O_y)@>>>\Cal F(k(y))
$
coincides with the map induced by the map
$
{\Cal O}_y @>\text{res}_y>> k(y))
$.
\enddefinition
\proclaim{4.6. Lemma-Definition}
Let
$
K
$
be a field containing the field
$
k
$
and let
$
K \subseteq L
$
be a finite field extension. Then the norm map
$
N_{L/K}:L^*@>>>K^*
$
takes the group
$
G_{L}(h)
$
into the group
$
G_{K}(h)
$.
Therefore the norm map
$
N_{L/K}
$
induces a map which we still denote by
$$
N_{L/K}:\Cal F(L)@>>>\Cal F(K).
$$
\endproclaim
\demo{Proof}
The Scharlau norm principle \cite{{\bf KMRT}, loc. cit.} states that
there is a natural inclusion
$
N_{L/K}(G_{L}(h))\subseteq G_{K}(h)
$,
whence the lemma.
\qed
\enddemo
\proclaim{4.7. Lemma}
Let
$
x
$
be a valuation of
$
K
$
and let
$
y
$
be a valuation of
$
L
$
extending
$
x
$.
Then the map
$
N_{\hat L_y/\hat K_x}:\Cal F(\hat L_y)@>>>\Cal F(\hat K_x)
$
takes
$
\Cal F(\hat \Cal O_y)
$
into
$
\Cal F(\hat \Cal O_x)
$.
\endproclaim
\demo{Proof}
The desired inclusion follows from the commutativity of the diagram
\comment
$$
\CD
\hat \Cal O_y^* @>>> {\longrightarrow} @>>>     {\longrightarrow} @>>> \hat L_y^*   \\
@VVV                 @.                         @.                     @VVV         \\
{\downarrow}@.       \Cal F(\hat \Cal O_y) @>>> \Cal F(\hat L_y)  @.   {\downarrow} \\
@VN_{y/x}VV          @.                         @VVNV                  @VVN_{y/x}V  \\
{\downarrow}@.       \Cal F(\hat \Cal O_x) @>>> \Cal F(\hat K_x)  @.   {\downarrow} \\
@VVV                 @.                         @.                     @VVV         \\
\hat \Cal O_x^* @>>> {\longrightarrow} @>>>     {\longrightarrow} @>>> \hat K_x^*,  \\
\endCD
$$
\endcomment
$$
\SelectTips{cm}{}
\xymatrix @C=3pc @R=1.5pc {
{\hat \Cal O_y^*} \ar"4,1"_{N_{y/x}} \ar@{->>}"2,2" \ar"1,4" & {} & {} & {\hat L_y^*} \ar"2,3" \ar "4,4"^{N_{y/x}} \\
{} & {\Cal F(\hat \Cal O_y)\ } \ar@{^{(}->}"2,3" & {\Cal F(\hat L_y)} \ar"3,3"^{N}  & {} \\
{} & {\Cal F(\hat \Cal O_x)\ } \ar@{^{(}->}"3,3" & {\Cal F(\hat K_x)} & {} \\
{\hat \Cal O_x^*} \ar"4,4" \ar"3,2"  & {} & {} & {\hat K_x^*\ ,} \ar"3,3"
}
$$
the surjectivity of the map
$
\hat \Cal O_y^* @>>> \Cal F(\hat \Cal O_y)$
and the injectivity of the map
$
\Cal F(\hat \Cal O_x)@>>>\Cal F(\hat K_x)
$
(see Remark 4.3).
\qed
\enddemo

\definition{4.8. Notation}
The map
$
\Cal F(\hat \Cal O_y)@>>>\Cal F(\hat \Cal O_x)
$
will be still denoted by
$
N_{y/x}
$.
\enddefinition
\definition{4.9. Notation}
Let
$
x
$
be a valuation of
$
K
$.
Set
$
\Cal F_x(L)=\underset{y/x}\to{\bigcap}\Cal F_y(L)
$.
\enddefinition
\proclaim{4.10. Lemma}
Let
$
x
$
be a valuation of
$
K
$.
Then
$
N_{L/K}(\Cal F_x(L))\subseteq \Cal F_x(K).
$
\endproclaim
\demo{Proof}
The desired inclusion follows from Lemma 4.8
and the commutativity of the diagram
\comment
$$
\CD
\underset{y/x}\to{\prod}\Cal F({\hat \Cal O_y}) @>>> {\longrightarrow} @>>>     {\longrightarrow} @>>> \underset{y/x}\to{\prod}\Cal F(\hat L_y)   \\
@VVV                 @.                         @.                     @VVV         \\
{\downarrow}@.       \Cal F_x(L) @>>> \Cal F(L)  @.   {\downarrow} \\
@V{\dsize\Pi}N_{y/x}VV          @.                         @VVN_{L/K}V                  @VV{\dsize\Pi}N_{y/x}V  \\
{\downarrow}@.       \Cal F_x(K) @>>> \Cal F(K)  @.   {\downarrow} \\
@VVV                 @.                         @.                     @VVV         \\
\Cal F(\hat \Cal O_x) @>>> {\longrightarrow} @>>>     {\longrightarrow} @>>> \Cal F(\hat K_x),  \\
\endCD
$$
\endcomment
$$
\SelectTips{cm}{}
\xymatrix @C=3pc @R=1.5pc {
{\underset{y/x}\to{\prod}\Cal F({\hat \Cal O_y})} \ar"4,1"_{{\dsize\Pi}N_{y/x}} \ar"1,4" & {} & {} & {\underset{y/x}\to{\prod}\Cal F(\hat L_y)} \ar"4,4"^{{\dsize\Pi}N_{y/x}} \\
{} & {\Cal F_x(L)\ } \ar@{^{(}->}"2,3" \ar"1,1" & {\Cal F(L)} \ar"3,3"^{N_{L/K}} \ar"1,4" & {} \\
{} & {\Cal F_x(K)\ } \ar@{^{(}->}"3,3" \ar"4,1" & {\Cal F(K)} \ar"4,4" & {} \\
{\Cal F(\hat \Cal O_x)} \ar"4,4" & {} & {} & {\Cal F(\hat K_x)\ ,}
}
$$
and the definition of
$
\Cal F_x(K)
$
(see 4.4).
\qed
\enddemo
\proclaim{4.11. Lemma}
Let
$
x
$
be a valuation of
$
K
$.
Then the diagram commutes.
$$
\CD
\Cal F_x(L) @>{\dsize\Pi}s_y>> \underset{y/x}\to{\prod}\Cal F(k(y))   \\
 @VN_{L/K}VV                      @VV{\dsize\Pi}N_{k(y)/k(x)}^{i(y/x)}V  \\
\Cal F_x(K) @>s_x>> \Cal F(k(x)),   \\
\endCD
$$
where
$
N_{k(y)/k(x)}:\Cal F(k(y)) @>>> \Cal F(k(x))
$
is the norm map for the field extension
$
k(y)/k(x)
$
and
$
N_{k(y)/k(x)}^{i(y/x)}
$
is its
$
i(y/x)
$-th power, where
$
i(y/x)
$
is the ramification index of
$
y
$
over
$
x
$, i.e.
$
i(y/x)=\text{ the length of }\Cal O_y/\eufm M_x\Cal O_y
$.
\endproclaim
\demo{Proof}
Consider the diagram
$$
\CD
\Cal F_x(L) @>>> \underset{y/x}\to{\prod}\Cal F(\hat \Cal O_y) @>{\dsize\Pi}\text{res}_y>> \underset{y/x}\to{\prod}\Cal F(k(y))   \\
 @VN_{L/K}VV  @VV{\dsize\Pi}N_{y/x}V @VV{\dsize\Pi}N_{k(y)/k(x)}^{i(y/x)}V  \\
\Cal F_x(K) @>i_*>> \Cal F(\hat \Cal O_x) @>\text{res}_x>> \Cal F(k(x))   \\
\endCD
$$
and observe that the left square commutes.
It remains to check that the right hand square commutes.
To do this it clearly suffices  to
check the commutativity of
$$
\CD
\Cal F(\hat \Cal O_y) @>\text{res}_y>> \Cal F(k(y))   \\
@VN_{y/x}VV @VVN_{k(y)/k(x)}^{i(y/x)}V  \\
\Cal F(\hat \Cal O_y) @>\text{res}_x>> \Cal F(k(x)).   \\
\endCD
$$
To see this we include it in a bigger one:
\comment
$$
\CD
\hat \Cal O_y^* @>>> {\longrightarrow} @>res_y>> {\longrightarrow} @>>> k(y)^*   \\
@VVV                 @.                            @.                   @VVV \\
{\downarrow} @. \Cal F(\hat \Cal O_y) @>res_y>> \Cal F(k(y))@. {\downarrow}   \\
@VN_{y/x}VV @VN_{y/x}VV @VVN_{k(y)/k(x)}^{i(y/x)}V @VVN_{k(y)/k(x)}^{i(y/x)}V  \\
{\downarrow} @.\Cal F(\hat \Cal O_y) @>res_x>> \Cal F(k(x)) @. {\downarrow}  \\
@VVV                 @.                            @.                   @VVV \\
\hat \Cal O_x^* @>>> {\longrightarrow} @>res_x>> {\longrightarrow} @>>> k(x)^* .  \\
\endCD
$$
\endcomment
$$
\SelectTips{cm}{}
\xymatrix @C=3pc @R=1.5pc {
{\hat \Cal O_y^*} \ar"1,4"^{\text{res}_y} \ar"4,1"_{N_{y/x}} \ar@{}"4,2"|{\txt{II}} \ar@{->>}"2,2"^{\rho} \ar@{}"2,4"|{\txt{I}\ \ \ } & {} & {} & {k(y)^*} \ar"2,3" \ar"4,4"^{N_{k(y)/k(x)}^{i(y/x)}} \\
{} & {\Cal F(\hat \Cal O_y)} \ar"2,3"^{\text{res}_y} \ar"3,2"_{N_{y/x}} \ar@{}"3,3"^{\txt{V}} & {\Cal F(k(y))} \ar"3,3"^{N_{k(y)/k(x)}^{i(y/x)}\ \ \txt{IV}} & {} \\
{} & {\Cal F(\hat \Cal O_y)} \ar"3,3"^{\text{res}_x} \ar@{}"4,3"^{\txt{III}} & {\Cal F(k(x))} & {} \\
{\hat \Cal O_x^*} \ar"4,4"^{\text{res}_x} \ar"3,2" & {} & {} & {k(x)^*.} \ar"3,3"
}
$$

The large square in this diagram commutes and squares I to IV commute as well
and the map
$
\rho
$
is surjective. Thus  square V commutes as well and the lemma is proved.
\qed
\enddemo
\proclaim{4.12. Proposition}
Let
$
K=k(t)
$
be the rational function field in one variable and
$
\Cal F_{un}(k(t))=\underset{x\in \Bbb A_k^1}\to{\bigcap}\Cal F_x(k(t))
$.
Then the canonical map
$$
\Cal F(k) @>>> \Cal F_{un}(k(t))
$$
is an isomorphism.
\endproclaim
\demo{Proof}
Injectivity is clear, because the composite
$
\Cal F(k) @>>> \Cal F_{un}(k(t)) @>s_0>> \Cal F(k)
$
coincides with the identity (here
$
s_0
$
is the specialization map at the point zero defined in 4.6).

It remains to check the surjectivity. Let
$
a\in \Cal F_{un}(k(t))
$.
Then by Lemma 4.5 the element
$
\partial (a)\in H_{et}^1(k(t),U_{A,\sigma})
$
is a class which for every
$
x\in \Bbb A_k^1
$
belongs to the image of
$
H_{et}^1(\Cal O_x,U_{A,\sigma})
$.
Thus by a lemma of Harder \cite{H},
$
\partial (a)
$
can be represented by an element
$
\xi \in H_{et}^1(k[t],U_{A,\sigma})
$,
where
$
k[t]
$
is the polynomial ring.
By  Harder's theorem \cite{H}, the map
$
H_{et}^1(k,U_{A,\sigma})@>>>H_{et}^1(k[t],U_{A,\sigma})
$
is an isomorphism. Then
$
\xi = \rho (\xi _0)
$
for an element
$
\xi _0\in H_{et}^1(k,U_{A,\sigma})
$.
Consider the diagram
\comment
$$
\CD
{}@.a@>>>\xi@>>>*@.{} \\
@.@.@.@.@. \\
1 @>>> \Cal F(k(t)) @>\partial>> H_{et}^1(k(t),U_{A,\sigma}) @>>> H_{et}^1(k(t),{\text{Sim}}_{A,\sigma}) @>>> 1 \\
@.     @AA\epsilon A               @AA\rho A              @AA\eta A            @. \\
1 @>>> \Cal F(k) @>\partial>> H_{et}^1(k,U_{A,\sigma}) @>>> H_{et}^1(k,{\text{Sim}}_{A,\sigma}) @>>> 1 , \\
@.@.@.@.@. \\
{}@.a_0@>>>\xi_0@.{}@.{} \\
\endCD
$$
\endcomment
$$
\SelectTips{cm}{}
\xymatrix @C=2pc @R=10pt {
{} & a \ar@{|->}[r] & {\xi} \ar@{|->}[r] & {\ast} & {} \\
1 \ar[r] & {\Cal F(k(t))} \ar[r]^{\partial\qquad} & {H_{et}^1(k(t),U_{A,\sigma})} \ar[r] & {H_{et}^1(k(t),{\text{Sim}}_{A,\sigma})} \ar[r] & 1 \\
{} & {} & {} & {} & {} \\
1 \ar[r] & {\Cal F(k)} \ar[r]^{\partial\qquad} \ar[uu]_{\epsilon} & {H_{et}^1(k,U_{A,\sigma})} \ar[r] \ar[uu]_{\rho} & {H_{et}^1(k,{\text{Sim}}_{A,\sigma})} \ar[r] \ar[uu]_{\eta} & 1 \\
{} & {a_0} \ar@{|->}[r] & {\xi_0 \ ,} & {} & {} }
$$
where all the mapping are canonical and all the vertical arrows
have trivial kernels. Since
$
\xi
$
goes to the trivial element in
$
H_{et}^1(k(t),{\text{Sim}}_{A,\sigma})
$,
one concludes that
$
\xi_0
$
goes to the trivial element in
$
H_{et}^1(k,{\text{Sim}}_{A,\sigma})
$.
Thus there exists an element
$
a_0 \in \Cal F(k)
$
such that
$
\partial (a_0)=\xi _0
$.
Clearly, one has
$
\epsilon(a_0)=a
$
(use the injectivity of the map
$
\Cal F(k(t)) @>>> H_{et}^1(k(t),U_{A,\sigma})
$
mentioned in Remark 4.3).
\qed
\enddemo

\proclaim{4.13. Theorem}
Let
$
L\supseteq K=k(t)
$
be a finite separable field extension and let
$
\Cal F_{un}(L)=\underset{y/x,\ x\in \Bbb A_k^1}\to{\bigcap}\Cal F_y(L)
$.
Then for an element
$
a\in \Cal F_{un}(L)
$
the following relation holds:
$$
\underset{y/0}\to{\dsize\Pi}N_{k(y)/k}(s_y(a)^{i(y/0)})=\underset{y/1}\to{\dsize\Pi}N_{k(y)/k}(s_y(a)^{i(y/1)}).
\tag{*}
$$
\endproclaim
\demo{Proof} By lemma 4.11, the element $ N_{L/K}(a) $ is in the
group $ \Cal F_{un}(K) $. Now by Lemma 4.12, the left hand side of
the relation (*) coincides with  $ s_0(N_{L/K}(a)) $,
where $ s_0:\Cal F_{un}(K)@>>>\Cal F(k) $ is the specialization
map (see Definition 4.6) at the point zero. The right hand side of
 (*) coincides with  $ s_1(N_{L/K}(a)) $,
where $ s_1 $ is the specialization map at  $1$. By
Proposition 4.13, there exists an element $ a_0\in \Cal F(k) $
whose image in $ \Cal F(k(t)) $ is equal to  $
N_{L/K}(a)\in \Cal F(k(t)) $. Thus
$$
s_0(N_{L/K}(a))=s_0(a_0)=a_0=s_1(a_0)=s_1(N_{L/K}(a)).
$$
The theorem is proved.
\qed
\enddemo

\proclaim{4.14. Corollary}
The Specialization Lemma (Theorem 3.3) holds.
\endproclaim
\demo{Proof}
We use notation of \S3.
Let
$
S\supseteq k[s]
$
be the integral extension of the polynomial ring in one variable and suppose
(as in the hypothesis of the Specialization Lemma) that
$
S
$
is an integral Dedekind domain. Let
$
L
$
be the quotient field of
$
S
$,
$
K=k(s)
$,
and
$
u\in S_{\eurm f}^*
$
for the element
$
\eurm f
$
from the hypotheses of the Specialization Lemma.

The element
$
{\overline u}\in \Cal F(L)
$
is $S$-unramified, i.e.
$
{\overline u}\in {\Cal F}_{un}(S)
$.
Thus
${\overline u}\in \Cal F_{un}(L)$.
Theorem 5.14 shows that the relation
$$
\underset{y/1}\to{\dsize\Pi}N_{k(y)/k}(s_y(\overline u)^{i(y/1)})=\underset{y/0}\to{\dsize\Pi}N_{k(y)/k}(s_y(\overline u)^{i(y/0)})
\tag{**}
$$
holds in $ \Cal F(k) $. It remains to check that the left hand
side of the relation (**) coincides with the element $ \overline
{N_{S_1/k}(u_1)} $ in $ \Cal F(k) $ and the right hand side of the
relation (**) coincides with the element $
\overline{N_{S_0/k}(u_0)}\cdot \overline{\epsilon(u)} $ in $ \Cal
F(k) $.

Let
$
S_{1,y}
$
be the localization at
$
y
$
of the Artinian ring
$
S_1=S/(s-1)S
$.
Clearly, the diagram
$$
\CD
S_{1,y}@>p_y>>k(y) \\
@AAA  @AAA \\
S_y@>>>\hat S_y , \\
\endCD
$$
(where all the mappings  are the  canonical ones) commutes.
For an element
$
v\in S_1
$
let
$
v_y
$
be its image in
$
S_{1,y}
$.
Now Lemma 4.5 and Definition 4.6 show that the element
$
p_y((u_1)_y)
$
coincides with the element
$
s_y(u)
$
in
$
\Cal F(k(y))
$.
Observe as well that
$
S_1=\underset{y/1}\to{\prod}S_{1,y}
$
and that the diagrams
$$
\CD
S_1^*@>\sim>>\underset{y/1}\to{\prod}S_{1,y}^* \\
@VN_{S_1/k}VV @VV{\dsize\Pi}N_{S_{1,y}/k}V \\
k^*@>\text{id}>> k^*, \\
\endCD
\qquad
\qquad
\qquad
\qquad
\CD
S_{1,y}^*@>p_y>>k(y)^* \\
@VN_{S_{1,y}/k}VV @VVN_{k(y)/k}^{i(y/1)}V \\
k^*@>\text{id}>> k^* . \\
\endCD
$$
commute. This proves the relation
$
\underset{y/1}\to{\dsize\Pi}N_{k(y)/k}(s_y(\overline u)^{i(y/1)})=\underset{y/1}\to{\dsize\Pi}N_{S_{1,y}/k}((u_1)_y)=N_{S_1/k}(u_1)
$
in
$
\Cal F(k)
$.
The relation
$
\underset{y/0}\to{\dsize\Pi}N_{k(y)/k}(s_y(\overline u)^{i(y/0)})=N_{S_0/k}(u_0)\cdot \epsilon(u)
$
in
$
\Cal F(k)
$
is proved similarly (use that
$
S/sS=k\times S_0
$
and the map
$
S@>>>k
$
is the augmentation
$
\epsilon:S@>>>k
$).
The Corollary is proved.
\qed
\enddemo

\head\S6.
Two lemmas
\endhead
Let $k$ be an infinite field and
$
\Cal O
$
 an essentially smooth  local
$
k
$-algebra.
\definition{5.1. Definition}
A {\it perfect triple\/}
$
(\Cal R {\underset{i}\to{\overset{\epsilon}\to{\rightleftarrows}}} \Cal O,f)
$
over
$
\Cal O
$
consists of a commutative
$
\Cal O
$-algebra
$
i:{\Cal O}@>>>{\Cal R}
$,
an augmentation map
$
\epsilon:{\Cal R}@>>>{\Cal O}
$
and an element
$
f\in {\Cal R}
$
which are subjected to the following conditions:
\roster
\item
$
\epsilon\circ i=id_{\Cal O}
$,
\item
the $ {\Cal O} $-algebra $ {\Cal R} $ is smooth at each prime $
\eufm p $ containing $ {\text{Ker}}(\epsilon) $,
\item
$ {\Cal R} $ is essentially $ k $-smooth and $ {\Cal R} $ is
domain,
\item
$
{\Cal R}/f{\Cal R}
$
is a finitely generated
$
{\Cal O}
$-module,
\item
there exists an element
$
t\in {\Cal R}
$
such that
$
{\Cal O}[t]
$
is the polynomial ring in one variable and
$
{\Cal R}
$
is a finitely generated
$
{\Cal O}[t]
$-module.
\endroster
\enddefinition
\remark
{5.2. Remark}
The condition (5) shows that
$
{\text{Spec}}\ \Cal R
$
is a relative curve over
$
{\text{Spec}}\ \Cal O
$.
The condition (3) shows that
$
\Cal R
$
is a regular ring. Since
$
\Cal O[t]
$
is a regular ring as well, a theorem of Grothendieck
\cite{{\bf Eis}, Corollary 18.17} together with (5) show that
$
\Cal R
$
is flat over
$
\Cal O[t]
$
and thus it is a finitely generated projective
$
\Cal O[t]
$-module.
\endremark

Let
$ (\Cal R{\underset{i}\to{\overset{\epsilon}\to{\rightleftarrows}}} \Cal O,f) $
a perfect triple. Let $ (A,\sigma) $ be an Azumaya
algebra with involution over $ \Cal R $ and let
$(A_0,\sigma_0)={\Cal O}\otimes_{\Cal R}(\Cal A,\sigma) $,
where $\Cal O $ is considered as an $ \Cal R $-algebra by means of the
augmentation
$ \epsilon $.
\proclaim{5.3. Lemma (Equating Lemma)}
There exist a quasi-finite \'etale extension
$ \widetilde j:\Cal R\hookrightarrow {\widetilde {\Cal R}} $
and a lifting
$\widetilde \epsilon :\widetilde {\Cal R} @>>> \Cal O $
of the augmentation
$ \epsilon $
(i.e.
$ \widetilde \epsilon \circ \widetilde j=\epsilon $)
and an isomorphism
$
\Phi: \widetilde
{\Cal R} \otimes_{\Cal R}(\Cal A,\sigma) @>>> \widetilde {\Cal R}
\otimes_{\Cal O}(A_0,\sigma_0)
$
of Azumaya algebras with
involutions over
$ \widetilde {\Cal R} $
such that the triple
$
(\widetilde {\Cal R} {\underset{\widetilde
i}\to{\overset{\widetilde {\epsilon}}\to{\rightleftarrows}}} \Cal O,\widetilde f)
$
with
$ \widetilde i=\widetilde j\circ i $
and
$
\widetilde f=\widetilde j(f)
$
is still perfect and the map
$ \Cal O \otimes_{\widetilde {\Cal R}}
\Phi: (A_0,\sigma_0)@>>>(A_0,\sigma_0)
$
is the identity.
\endproclaim

\demo{Proof}
The required
quasi-finite \'etale extension
$ \widetilde j:\Cal R\hookrightarrow {\widetilde {\Cal R}} $
and lifting
$\widetilde \epsilon :\widetilde {\Cal R} @>>> \Cal O $
of the augmentation
$ \epsilon $
and the isomorphism
$
\Phi: \widetilde
{\Cal R} \otimes_{\Cal R}(\Cal A,\sigma) @>>> \widetilde {\Cal R}
\otimes_{\Cal O}(A_0,\sigma_0)
$
of Azumaya algebras with
involutions
are constucted using geometric terminology in
\cite{{\bf Oj-P1}, Proof of  8.1}.

To see this set
$\Cal X={\roman{Spec}}(\Cal R)$, $U={\roman{Spec}}(\Cal O)$
and consider the
morphisms $p: \Cal X \to U$ and $\Delta: U \to \Cal X$ induced by the ring
homomorphisms $i$ and $\epsilon$. Let
$q: \Cal X \to U \times {\Bbb A}^1$
be the finite surjective $U$-morphism
corresponding to the integral extension
$\Cal O [t] \subset \Cal R$. Let
$\Cal Z \subset \Cal X$
be the vanishing locus of $f$.

Now consider certain scheme morphisms from
\cite{{\bf Oj-P1}, Proof of  8.1}. Namely, consider
the quasi-finite \'etale morphism
$\tilde \Cal X \to \Cal X$
which is the composition of the finite surjective etale morphism
$\pi: \tilde \Cal X \to \Cal W$
and the open inclusion
$\Cal W \subset \Cal X$.
Consider the section
$\widetilde \Delta: U \to \tilde \Cal X$
and the isomorphism of Azumaya algebras with involutions
$\Phi$ from
\cite{{\bf Oj-P1}, Proof of  8.1}.
Recall that
$\widetilde \Delta^*(\Phi)$
is the identity,
$\Delta = \pi \circ \widetilde {\Delta}$,
$\Delta (U) \subset \Cal W$,
$\Cal Z \subset \Cal W$,
and that there is a finite surjective $U$-morphism
$r: \Cal W \to U \times {\Bbb A}^1$.

Let
$\widetilde j : \Cal R \hookrightarrow \widetilde {\Cal R}$
be the inclusion induced by
$\tilde \Cal X \to \Cal X$
and
$\widetilde \epsilon: \widetilde {\Cal R} \to \Cal O$
the $\Cal O$-augmentation induced by
$\widetilde \Delta: U \to \tilde \Cal X$.
We claim that
$\widetilde j$,
$\widetilde \epsilon$
and $\Phi$
satisfy the Lemma.

In fact,
$\Cal O \otimes_{\widetilde {\Cal R}} \Phi = \widetilde \Delta^*(\Phi)$
is the identity.
The relation
$ \widetilde \epsilon \circ \widetilde j=\epsilon $
follows from the equality
$\Delta = \pi \circ \widetilde {\Delta}$
mentioned just above.
It remains to check that the triple
$
(\widetilde {\Cal R} {\underset{\widetilde
i}\to{\overset{\widetilde {\epsilon}}\to{\rightleftarrows}}} \Cal O,\widetilde f)
$
is perfect.
To check this note that
$\widetilde {\epsilon} \circ \widetilde {i} =
\widetilde {\epsilon} \circ \widetilde {j} \circ i =
\epsilon \circ i = id_{\Cal O}$.
The $\Cal O$-algebra $\widetilde {\Cal R}$ is smooth at
each prime containing
$\text{Ker} (\widetilde {\Delta})$
because
$\Delta= \pi \circ \widetilde {\Delta}$
(with $\pi$ an \'etale morphism ) and
$p: \Cal X \to U$
is smooth along $\Delta (U)$.
The $k$-algebra $\widetilde {\Cal R}$ is essentially smooth
because the $k$-algebra $\Cal R$ is essentially smooth and
$\widetilde j: \Cal R \to \widetilde {\Cal R}$
is etale.
The vanishing locus
$\widetilde {\Cal Z} \subset \widetilde {\Cal X}$ of $\widetilde f$
is finite over $U$
because
$\Cal Z \subset \Cal W$
and
$\pi: \tilde \Cal X \to \Cal W$
is finite. Since $\widetilde {\Cal Z}$ is finite over $U$
the $\Cal O$-module ${\Cal R}/f{\Cal R}$ is finitely generated.
It remains to check that there is a finite surjective $U$-morphism
$\widetilde {\Cal X} \to U \times {\Bbb A}^1$. For that consider the finite
surjective morphism
$r: \Cal W \to U \times {\Bbb A}^1$
and take the composition
$r \circ \pi: \widetilde {\Cal X} \to U \times {\Bbb A}^1$.
\qed
\enddemo

Let $ R $ be a commutative $ k $-algebra and let $ (A,\sigma) $ be
an Azumaya algebra with involution over $ R $. Let $
G={\text{Sim}}_{A,\sigma} $ be the similitude group of $
(A,\sigma) $ and let $ \mu:G @>>> \Bbb G_m $ be a group
homomorphism which takes a similitude $ \alpha $ to its similarity
factor $ \mu(\alpha)={\alpha}^{\sigma}\cdot \alpha $. Observe that
$ \mu(G(S))=G_S(h) $ from 3.1. \qed

\definition{5.4. Notation}
For every commutative
$
R
$-algebra
$
S
$
denote by
$
{\Cal F}(S)
$
the group
$
S^*/\mu(G(S))
$.
An
$
R
$-algebra homomorphism
$
S@>\alpha>>T
$
clearly induces a group map
$
{\Cal F}(S)@>\alpha_*>>{\Cal F}(T)
$.
For an element
$
u\in S^*
$
we shall write
$
\overline u
$
for its image in
$
{\Cal F}(S)
$.
The homomorphism
$
\alpha_*
$
takes
$
\overline u
$
to
$
\overline {\alpha(u)}
$.
\enddefinition
\definition{5.5. Definition}
Let
$
S
$
be an
$
\Cal R
$-algebra which is a domain with the quotient field
$
K
$
and let
$
\eufm p
$
be a height $1$ prime ideal in
$
S
$.
An element
$
{\eurm v}\in {\Cal F}(K)
$
is called unramified at
$
\eufm p
$
iff
$
\eurm v
$
belongs to the image of
${\Cal F}(S_{\eufm p})$
in
${\Cal F}(K)$.
An element
${\eurm v}\in {\Cal F}(K)$
is called
$S$-unramified if it is unramified at each height 1 prime
$\eufm p$
in
$S$.
\enddefinition
\proclaim{5.6. Lemma (Unramifiedness Lemma)}
Let $R$ and $S$ be domains with  quotient fields $K$ and $L$
respectively.
Let
$
R@>\alpha>>S
$
be an injective flat homo\-morph\-ism of finite type and
let
$
\beta :K@>>>L
$
be the induced inclusion of the quotient fields.
Then for each localization $T \supset S$ of $S$
the map
$
\beta _*:{\Cal F}(K)@>>>{\Cal F}(L)
$
takes
$
S
$-unramified elements to
$
T
$-unramified elements.
\endproclaim
\demo{Proof}
Let
$\eurm v\in K^*$
and let
$\eufm r$
be  height $1$ primes of $T$. Then
$\eufm q = S \cap \eufm r$
is a height $1$ prime of $S$. Let
$\eufm p = R \cap \eufm q$.
Since the $R$-algebra $S$ is flat of finite type one has
$\text{ht}(\eufm q) \geq \text{ht}(\eufm p)$. Thus
$\text{ht}(\eufm p)$
is $1$ or $0$. The commutative diagram
$$
\CD
\Cal F(K) @>>> \Cal F(L) \\
@AAA @AAA \\
\Cal F(R_{\eufm p}) @>>> \Cal F(T_{\eufm r}) \\
\endCD
$$
shows that the class
$\overline {\beta (\eurm v)}$
is in the image of
$\Cal F(T_{\eufm r})$.
Whence the class
$\overline {\beta (\eurm v)} \in \Cal F(L)$
is $T$-unramified. The Lemma follows.
\qed
\enddemo
\head\S6.
Relative Specialization Lemma
\endhead
Let $\Cal O$ be a regular local ring containing an infinite field
$k$ and which is an essentially smooth $k$-algebra. Let $K$ be the
quotient field of $\Cal O$. Let $ (\Cal R
{\underset{i}\to{\overset{\epsilon}\to{\rightleftarrows}}} \Cal O,
f) $ be a perfect triple. Denote by $\epsilon_K: \Cal R_K= \Cal R
\otimes_{\Cal O}  K \to K$ the homomorphism $\epsilon \otimes_{\Cal O} K$.
We will consider $\Cal O$ and $K$ as  $\Cal R$-algebras via
$\epsilon$ and $\epsilon_K$ respectively. So for an Azumaya
algebra with involution $(\Cal A, \sigma)$ over $\Cal R$ it makes
sense to speak about the groups $\Cal F(\Cal O)$ and $\Cal F(K)$
(see Definition 5.5).

\proclaim{6.1. Lemma (Relative Specialization Lemma)}
Let
$
(\Cal R {\underset{i}\to{\overset{\epsilon}\to{\rightleftarrows}}} \Cal O, f)
$
be a perfect triple and
$(\Cal A,\sigma)$
an Azumaya algebra with involution over $\Cal R$.
Let $\Cal K$ be the quotient field of $\Cal R$ and
let $u \in \Cal R_f^*$ be a unit such that the class
$\bar u \in \Cal F(\Cal K)$
is
$\Cal R$-unramified. If
$\epsilon (f)\not= 0$
then the class
$\overline {\epsilon_K (u \otimes 1)} \in \Cal F(K)$
can be lifted to
$\Cal F(\Cal O)$.
\endproclaim

\demo{Proof} Set
$ (A_0, \sigma_0)=(\Cal O \otimes_{\Cal R} \Cal A, \Cal O \otimes_{\Cal R} \sigma) $,
where $ \Cal O $ is an $\Cal R $-algebra by means of $ \epsilon $.
Set
$ ({\Cal A}_0,{\sigma}_0)=(\Cal R \otimes_{\Cal O} A_0, \Cal R \otimes_{\Cal O} \sigma_0) $,
where $ \Cal R $ is regarded  as an
$\Cal O$-algebra by means of the map $i$.
There are two Azumaya
algebras with involutions
$(\Cal A, \sigma)$ and
$({\Cal A}_0,{\sigma}_0)$
over $\Cal R$.
Their scalar extensions
$ (\Cal A,\sigma)\otimes_{\Cal R} \Cal O $ and
$ ({\Cal A}_0,{\sigma}_0)\otimes_{\Cal R} \Cal O $
tautologically coincides
because the composite map
$ \Cal O @>i>> \Cal R @>\epsilon>> \Cal O $
is the identity. Thus by the Equating Lemma 5.3 one can find a
quasi-finite \'etale extension
$ j:\Cal R\hookrightarrow {\widetilde {\Cal R}} $
and a lifting
$ \widetilde \epsilon: \widetilde {\Cal R} @>>> \Cal O $
of the augmentation
$ \epsilon$ and an isomorphism
$ \Phi: (\widetilde {\Cal A},\widetilde {\sigma}) @>>> (\widetilde {\Cal A_0},\widetilde {\sigma_0}) $
of Azumaya algebras with involutions over
$ \widetilde {\Cal R} $
such that
$ (\widetilde {\Cal R} {\underset{\widetilde i}\to
{\overset{\widetilde {\epsilon}}\to {\rightleftarrows}}} \Cal O, {j}(f)) $
is still a perfect triple and the  isomorphism
$ \Cal O \otimes_{\widetilde {\Cal R}} \Phi $
is the identity. Here
$(\widetilde {\Cal A},\widetilde {\sigma})= \widetilde {\Cal R} \otimes_{\Cal R}(\Cal A,\sigma) $,
$ (\widetilde {\Cal A_0},\widetilde {\sigma_0})= \widetilde {\Cal R} \otimes_{\Cal O}(A_0,\sigma_0) $
and
$ \Cal O $
is regarded as an
$ \widetilde {\Cal R} $-algebra
by means of
$ \widetilde \epsilon :\widetilde {\Cal R} @>>> \Cal O $.

Denote by
$\tilde \epsilon_K: \tilde \Cal R_K= \Cal R \otimes_{\Cal O}  K \to K$
the augmentation
$\tilde \epsilon \otimes_{\Cal O} K$.
Set
$\tilde f = j(f)$
and
$\tilde u = j(u) \in \tilde \Cal R_{\tilde f}^*$.
Since
$\tilde \epsilon_K (\tilde u \otimes 1)=\epsilon_K (u \otimes 1)$
it suffices to check that the class
$\overline {\tilde \epsilon_K (\tilde u \otimes 1)}$
can be lifted to $\Cal F(\Cal O)$.

Let $\tilde {\Cal K}$ be the quotient field of $\tilde \Cal R$. By
the Unramifiedness Lemma the class $\overline {\tilde u} \in \Cal
F(\tilde {\Cal K})$ is $\tilde \Cal R$-unramified. So replacing $
(\Cal R {\underset{i}\to{\overset{\epsilon}\to{\rightleftarrows}}}
\Cal O, f) $, $(\Cal A, \sigma)$ and $u$ by $ (\widetilde {\Cal R}
{\underset{\widetilde i}\to {\overset{\widetilde
{\epsilon}}\to{\rightleftarrows}}} \Cal O, {j}(f)) $, $(\tilde
\Cal A, \tilde \sigma)$ and $\tilde u$ we may assume that $ (\Cal
R {\underset{i}\to{\overset{\epsilon}\to{\rightleftarrows}}} \Cal
O, f) $, is a perfect triple, $(\Cal A, \sigma)= \Cal R
\otimes_{\Cal O} (A_0, \sigma_0)$ for an Azumaya algebra with
involution $(A_0, \sigma_0)$ over $\Cal O$, and $u \in \Cal R_f^*$
is such that the class $\bar u \in \Cal F(\Cal K)$ is $\Cal
R$-unramified. We must check that the class $\overline
{\epsilon_K(u_K)} \in \Cal F(K)$ can be lifted in $\Cal F(\Cal
O)$.

Since the triple
$
(\Cal R {\underset{i}\to{\overset{\epsilon}\to{\rightleftarrows}}} \Cal O, f)
$
is perfect, the geometric presentation lemma
\cite{{\bf Oj-P}, Lemma 5.2}
shows that
one can choose an element
$
s\in \Cal {R}
$
such that the extension
$
\Cal {R}\supseteq \Cal {O}[s]
$
is finite, the ring
$
\Cal O[s]
$
is the polynomial ring in one variable over
$
\Cal O
$
and the following holds:

\roster
\item
$
(1 - s)\Cal {R} + {f}\Cal {R}=\Cal {R}
$,
\item
$
s\Cal {R}={\text{Ker}}({\epsilon}){\cap} J
$
for a certain ideal
$
J
$
and
\item
$
J+ f\Cal {R}=\Cal {R}
$
and
\item
the map
$
\Cal {R}/s \Cal {R}\longrightarrow \Cal {R}/{\text{Ker}}({\epsilon}) \times {{\Cal {R}}/{J}}=
\Cal O \times \Cal {R}/J
$
is an isomorphism.
\endroster

Since
$
\Cal {R}
$
and
$
\Cal {O}[s]
$
are both essentially smooth
$
k
$-algebras (and thus regular rings) and since the extesion
$
\Cal {R}
$
over
$
{\Cal O}[s]
$
is finite, a theorem of Grothendieck
\cite{{\bf Eis}, Corollary 18.17}
shows that
$
\Cal {R}
$
is a flat
$
{\Cal O}[s]
$-module. Therefore
$
\Cal {R}
$
is a finitely generated projective
$
{\Cal O}
$-module.
Thus
$
{\Cal {R}}_1 = {\Cal {R}} / (1-s){\Cal {R}}
$
and
$
{\Cal {R}}_0 =\Cal {R}/{J}
$
are finitely generated projective
${\Cal O}$-modules.

Consider the elements
$u_1=u {\text{ mod }} (1-s){\Cal {R}}_{f}$
in
${\Cal {R}}_{1, {f}}^*$
and
$u_0=u {\text{ mod }} {J}_f$
in
${\Cal {R}}_{0, {f}}^*$.
By (1) and (3) one has
$
{\Cal {R}}_i={\Cal {R}}_{i, f}
$
and thus
$
u_i\in {\Cal {R}}_i^*
$
$
(i=0,1)
$.
Since
$
{\Cal {R}}_1
$
and
$
{\Cal {R}}_0
$
are finitely generated projective
$
{\Cal O}
$-modules, there are the norm mappings
$
N_{{\Cal {R}}_i/{\Cal O}}:{\Cal {R}}_i^*\longrightarrow{\Cal O}^*
$
$
(i=0, 1)
$
given by
$
(v\mapsto {\text{det}}(\text{mult. by }v))
$.
Set
$$
\phi(u)=N_{{\Cal {R}}_1/{\Cal O}}(u_1)\cdot {N_{{\Cal {R}}_0/{\Cal O}}(u_0^{-1})}\in {\Cal O}^*
\subseteq K^* \text{ .}
$$
\proclaim{Claim}
$
\overline {\phi(u)}=\overline {\epsilon_K (u_K)}
$
in the group
$
K^*/\mu(G(K))={\Cal F}(K)
$.
\endproclaim
Since
$
\phi(u)\in{\Cal O}^*
$,
the Claim clearly completes the proof of purity.
The rest of the section is devoted to the proof of the Claim.

Set
$\Cal R_K= K \otimes_{\Cal O} \Cal R$
and
$u_K=1 \otimes u \in \Cal R_K$.
Set
$\Cal R_{i,K}= K \otimes_{\Cal O} \Cal R_i$
and
$u_{i,K}=1 \otimes u_i \in \Cal R_K \in \Cal R_{i, K}^*$.
Finally set
$\Cal R_{f, K} = \Cal R_{K, 1 \otimes f}$.
Clearly it suffices to prove the relation
$$
{\overline {{\epsilon}_K(1\otimes u)}}={\overline {N_{\Cal R_{1, K}/K}(1\otimes u_1)}}\cdot
{\overline {N_{\Cal R_{0, K}/K}(1\otimes u_0)^{-1}}}
\tag{\dag}
$$
in the group
$
\Cal {F}(K)
$.
The relation (\dag) will be checked below in this proof applying the
Specialization Lemma (Theorem 3.3) to the integral extension
$
\Cal R_K \supseteq K\otimes_{\Cal O}{\Cal O}[s]=K[s]
$
and the Azumaya algebra with involution
$(A_0, \sigma_0) \otimes_{\Cal O} K$
over $K$.

Check the hypotheses of the Specialization Lemma.
Since $\Cal R$ is regular domain and $\Cal R_K$
is its localization $\Cal R_K$ is a regular domain as well.
Since
$\Cal R_K$
is an integral extension of the polynomial ring
$K[s]$,
the dimension of
$\Cal R$
is one. A regular domain of dimension $1$ is a Dedekind domain.
\underbar {Thus}
$\Cal R$
\underbar {is a Dedekind domain.}

The class
$\bar u \in \Cal F(\Cal K)$
of the element
$u \in \Cal R_f^*$
is $\Cal R$-unramified. Thus the class
$\bar u_K \in \Cal F(\Cal K)$
of the element
$u_K \in \Cal R_{f, K}^*$
is $\Cal R_K$-unramified.

Now check that the element
$1 \otimes f \in \Cal R_K$
is coprime with both $s$ and $s-1$ in
$\Cal R_K$. Recall the conditions $(1)$ to $(4)$ mentioned above in this proof.
The element $1 \otimes f$ is coprime with $(s-1)$
by  condition $(1)$.
The element
$\epsilon_K (1 \otimes f) = \epsilon (f)_K$
is non-zero in $K$ by the very assumption on $f$.
The element $1 \otimes f$ is coprime with the ideal $J_K$
by  condition $(3)$.
Thus $1 \otimes f$ is coprime with $s$ by condition $(4)$.
We already checked that the class $\bar u$ is $\Cal R_K$-unramified.
Thus by Theorem $3.3$ the relation $(\dag)$ holds in $\Cal F(K)$.
The Claim is proved. The Relative Specialization Lemma follows.

\enddemo
\head\S7.
Geometric case of the Purity Theorem
\endhead
Under the notation of 5.4 and 5.5 the following theorem holds.
\proclaim{7.1. Theorem}
Let $\Cal O$ be a local, essentially smooth algebra over a field $k$
and let $K$ be its quotient field. Let $(A, \sigma)$ be an Azumaya algebra
with involution over $\Cal O$ and let $\eurm v \in K^*$ be such that the class
$\bar \eurm v \in \Cal F(K)$
is $\Cal O$-unramified. Then $\bar \eurm v$ can be lifted in $\Cal F(\Cal O)$.
\endproclaim
\demo{Proof}
We begin with the case of an infinite field $k$.
By assumption there exist a smooth $d$-dimensional $k$-algebra
$R=k[t_1, \dots, t_n]$
and a prime ideal $\eufm p$ of $R$ such that $A=R_{\eufm p}$.
We first reduce the proof to the case in which $\eufm p$ is maximal.
To do this we choose a maximal ideal $\eufm m$ containing $\eufm p$.
Since $k$ is infinite, by a standart general position argument we can find $d$
algebraically independent elements $X_1, X_2, \dots, X_d$ such that
$R$ is finite over $k[X_1, \dots, X_d]$ and \'etale at $\eufm m$.
After a linear change of coordinates we may assume that
$R/\eufm p$ is finite over
$B=k[X_1, \dots, X_m]$,
where $m$ is the dimension of $R/\eufm p$.
Clearly $R$ is smooth over $B$ at $\eufm m$ and thus,
for some $h \in R - \eufm m$, the localization $R_h$ is
smooth over $B$.
Let $S$ be the set of nonzero elements of $B$,
$k^{\prime}=S^{-1}B$ the field of fractions of $B$ and
$R^{\prime}=S^{-1}R_h$.
The prime ideal
$\eufm p^{\prime}=S^{-1}\eufm p_h$
is maximal in $R^{\prime}$, the
$k^{\prime}$-algebra $R^{\prime}$ is smooth and
$A=R^{\prime}_{\eufm p^{\prime}}$.

From now on and till the end of the proof of Theorem 8.1
we assume that
$\Cal O = \Cal O_{X,x}$
is the local ring of a closed point $x$ of a smooth
$d$-dimensional irreducible affine variety $X$ over $k$.

Replacing $X$ by a sufficiently small affine neighbourhood of $x$
we may assume that
\roster
\item
the algebra with involution
$(A, \sigma)$
is defined over $k[X]$ and is an Azumaya algebra with involution
already over $k[X]$,
\item
the element $\eurm v$ is a unit in $k[X]_g$ for certain nonzero element
$g \in k[X]$,
\item
the class $\bar \eurm v \in \Cal F(K)$ is $k[X]$-unramified.
\endroster
We must prove that $\bar \eurm v$ can be lifted in $\Cal F(\Cal O)$.

By Quillen's trick there exists a polinomial subalgebra
$k[t_1,t_2, \dots , t_n]$ in $k[X]$
such that the algebra $R=k[X]$ is finite over
$k[t_1,t_2, \dots , t_n]$,
the algebra $R$ is smooth over
$k[t_1,t_2, \dots , t_{n-1}]$
at the maximal ideal
$\eufm m$
and the
$k[t_1,t_2, \dots, t_{n-1}]$-module $R/fR$
is finite. Set
$P=k[t_1,t_2, \dots , t_{n-1}]$,
$\Cal {R}=\Cal {O} \otimes_{P} R$,
consider ring homomorphisms
$j: R \to \Cal {R}$,
$i: \Cal {O} \to \Cal{R}$
and
$\epsilon: \Cal{R} \to \Cal {O}$
given by
$j(a)=1 \otimes a$,
$i(b)=b \otimes 1$
and
$\epsilon (a \otimes b)=ab$
respectively.

We claim that
$
(\Cal R {\underset{i}\to{\overset{\epsilon}\to{\rightleftarrows}}} \Cal O, f)
$
with
$f=j({\eurm f})$
is a perfect triple (see 6.1 for definition).
This is checked in
[Oj-P]
using  geometric terminology.
This perfect triple fits in the diagram

\comment
$$
\CD
\Cal {R} @<j<< R \\
\epsilon \downarrow @AAiA \swarrow \\
\Cal {O}   @.{} \\
\endCD
$$
\endcomment
$$
\SelectTips{cm}{}
\xymatrix @C=3pc @R=2pc {
{\Cal R} \ar@<-2pt>[d]_{\epsilon} & R \ar[dl]^{can} \ar[l]_j \\
{\Cal O} \ar@<-2pt>[u]_i
}
$$
with the localization map $can$. Clearly $can=\epsilon \circ j$.

Set $(\Cal A, \sigma)=\Cal R \otimes_{k[X]} (A, \sigma)$ and
$u=j(\eurm v) \in \Cal R_f^*$.
Let $\Cal K$ be the quotient field of $\Cal R$.
By the Unramifiedness Lemma the class
$\bar u \in \Cal F(\Cal K)$
is
$\Cal R$-unramified.
Since $\epsilon (f) = \epsilon(j(f))=f$ is nonzero element of $\Cal O$
we are under the hypotheses of the Relative Specialization Lemma.
Thus the class
$\overline \epsilon_K (u) \in \Cal F(K)$
can be lifted in $\Cal F(\Cal O)$. It remain to note that
$$
\epsilon_K (u) = \epsilon_K (j(\eurm v)) =\eurm v \in K.
$$
Thus the class
$\bar \eurm v \in \Cal F(K)$
can be lifted in
$\Cal F(\Cal O)$.

Now suppose that $k$ is finite. So $\Cal O$ is a local essentially
smooth $k$-algebra with maximal ideal $\eufm m$. Let
$\eurm v \in K^*$
be such that the class
$\bar \eurm v \in \Cal F(K)$
is $\Cal O$-unramified.
Let $p^m$ be the cardinality of the algebraic closure of $k$ in
$A/\eufm m$
and $s$ be an odd integer greater than $2$ and prime to $m$.
For any $i$ let $l_i$ be the field
(in some fixed algebraic closure of $k$)
of degree $s^i$ over $k$. Let $l$ be the union of all $l_i$. Since
$l \otimes_k (\Cal O/\eufm m)$
is still a field,
$R=l \otimes_k \Cal O$
is a local essentially smooth algebra over the infinite field $l$. Let
$L= l \otimes_k K$
be its field of fractions. The image
$\bar \eurm v_{L}$ of
$\bar \eurm v$ in
$\Cal F(L)$ is
$R$-unramified.
In fact, let
$\eufm q$
be a hight-one prime of $R$ and
$\eufm p =\Cal O \cap \eufm q$.
By assumption
$\bar \eurm v$
is in the image of
$\Cal F(\Cal O_{\eufm p})$
and since
$\Cal O_{\eufm p} \to L$
factors through
$R_{\eufm q}$
the class
$\bar \eurm v_L$
is in
$\Cal F(R_{\eufm q})$
for every
$\eufm q$.
We can now find a finite subfield $l^{\prime}$ of $l$, and for
$\Cal O^{\prime}=l^{\prime} \otimes_k \Cal O$,
a $\eurm v^{\prime} \in \Cal O^{\prime}$
which maps to
$\bar \eurm v_L$.
Let $K^{\prime}$ be the field of fractions of
$\Cal O^{\prime}$.
Further enlarging $l^{\prime}$ we may assume that the images
$\bar \eurm v$ and $\bar \eurm v^{\prime}$ in
$\Cal F(K^{\prime})$ coincide. Consider the diagram
$$
\CD
\Cal O^* @>>> \Cal (O^{\prime})^* @>{N}>> \Cal O^* \\
@VVV @VVV @VV{\alpha}V \\
\Cal F(K) @>>> \Cal F(K^{\prime}) @>{\bar N}>> \Cal F(K) . \\
\endCD
$$
where $\bar N$ is the norm map
(it is well-defined by the Scharlau norm principle).
Since the composition of the horizontal maps is the identity,
we have
$\alpha \circ N(\eurm v^{\prime})=\bar \eurm v$
in $\Cal F(K)$. Thus $\bar \eurm v$ is indeed in the image of
$\Cal O^*$. Theorem 7.1 is proved.

\enddemo

\head\S8.
Proof of the Purity Theorem
\endhead
\demo{Proof of Theorem 1.3}
Let $k$ be the prime subfield of the ring $R$.
By Popescu's theorem [P], [Sw]
$R=\raise-6pt\hbox{$\buildrel{\displaystyle\lim}\over{\to}$}
R_\alpha$
(a filtered direct limit),
where $R_{\alpha}$'s are smooth $k$-algebras.
We first observe that we may replace the direct system
of the $R_{\alpha}$'s by a system of essentially smooth
local $k$-algebras. In fact, if $\eufm m$ is the maximal ideal
of $R$, we can replace each
$R_{\alpha}$ by
$(R_{\alpha})_{\eufm p_{\alpha}}$
where
$\eufm p_{\alpha} = \eufm m \cap R_{\alpha}$.
Note that in this case the canonical morphisms
$\phi_{\alpha}: R_{\alpha} \to R$
are local and that every
$R_{\alpha}$
is a regular local ring thus in particular a factorial ring.

Now let $K$ be the field of fractions of $R$ and, for each $\alpha$,
let $K_{\alpha}$ be the field of fractions of $R_{\alpha}$. The ideal
$\eufm r_{\alpha} = \ker (\phi_{\alpha})$
is prime. Set
$S_{\alpha}=(R_{\alpha})_{\eufm r_{\alpha}}$.
Note that the
$S_{\alpha}$'s
form a direct system of regular local rings with
$K=\raise-6pt\hbox{$\buildrel{\displaystyle\lim}\over{\to}$}
S_\alpha$
(a filtered direct limit).

We may assume that there exists an index $\alpha$ and
an Azumaya algebra with involution
$(A_{\alpha}, \sigma_{\alpha})$
over $R_{\alpha}$ such that
$(A, \sigma) = (A_{\alpha}, \sigma_{\alpha}) \otimes_{R_{\alpha}} R$.
Replacing the direct system of indeces $\alpha$'s by
the subsystem of indeces $\beta$ satisfying
$\beta \geq \alpha$
we may assume that we are given with a direct system of Azumaya algebras
with involutions
$(A_{\alpha}, \sigma_{\alpha})$
over the $R_{\alpha}$'s such that
$(A, \sigma)=\raise-6pt\hbox{$\buildrel{\displaystyle\lim}\over{\to}$}
(A_\alpha, \sigma_\alpha)$.

Let
$G_{\alpha}=\text{Sim}_{A_{\alpha},\sigma_{\alpha}}$.
Then one has
$G(R)=\raise-6pt\hbox{$\buildrel{\displaystyle\lim}\over{\to}$}
G_{\alpha}(R_{\alpha})$
and
$G(K)=\raise-6pt\hbox{$\buildrel{\displaystyle\lim}\over{\to}$}
G_{\alpha}(S_{\alpha})$.
Let
$\mu_{\alpha}: G_{\alpha} \to \Bbb G_m$
be the group morphism which takes a similitude
to its similarity factor (see the Introduction).

Let $\bar a \in \Cal F(K)$ be an $R$-unramified class.
We may represent $\bar a$ by a unit $a \in R^*_f$,
where $0\not= f \in R$. Let
$f= p_1p_2 \dots p_n$
be a prime decomposition of $f$ in $R$.
Since $\bar a$ is $R$-unramified for every index
$i=1,2, \dots ,n$
there exist elements
$h_i \in R - p_iR$
and
$a_i \in R^*_{h_i}$
and
$g_i \in G(K)$
such that
$a= a_i \mu (g_i)$.

We can now choose an index $\alpha$, elements
$p_{\alpha, i}$ and $h_{\alpha, i}$
with
$\phi_{\alpha}(p_{\alpha, i})=p_i$
and
$\phi_{\alpha}(h_{\alpha, i})=h_i$.
Set
$f_{\alpha}=p_{\alpha, 1}p_{\alpha, 2} \dots p_{\alpha, n}$.
Since
$\phi_{\alpha}(f_{\alpha}) = f \not= 0$
and
$\phi_{\alpha}(h_{\alpha, i}) = h_i \not= 0$
one has the inclusions
$R_{\alpha, f_{\alpha}} \subset S_{\alpha}$
and
$R_{\alpha, h_{\alpha, i}} \subset S_{\alpha}$.
Futher enlarging the index $\alpha$ we can choose
$a_{\alpha} \in R^*_{\alpha, f_{\alpha}}$
and elements
$a_{\alpha, i} \in R^*_{\alpha, h_{\alpha, i}}$
and
$g_{\alpha, i} \in G_{\alpha}(S_{\alpha})$
which are preimages of the $a$ and
the $a_i$'s and the $g_i$'s respectively. Having choosen
these preimages consider the relations
$$
a_{\alpha}=a_{\alpha, i}\mu_{\alpha}(g_{\alpha, i})
$$
in $S_{\alpha}$. Since they hold over $K$,
we may assume, after replacing $\alpha$ by
some larger index, that they hold over $S_{\alpha}$.
We claim that the class
$\bar a_{\alpha} \in \Cal F(K_{\alpha})$
is
$R_{\alpha}$-unramified.

To prove this  note first that each
$p_{\alpha, i}$
is prime. In fact,
$\phi_{\alpha}(p_{\alpha, i})=p_i$,
the element $p_i$ is prime and
$\phi_{\alpha}$
is a local homomorphism of local factorial rings.
Thus $p_{\alpha, i}$ is  indeed prime. Since
$a_{\alpha} \in R^*_{\alpha, f_{\alpha}}$
and
$f_{\alpha}=p_{\alpha, 1}p_{\alpha, 2} \dots p_{\alpha, n}$
the class
$\bar a_{\alpha}$
can be ramified at most at one of
the $p_{\alpha, i}$'s.
However the relations
$a_{\alpha}=a_{\alpha, i}\mu_{\alpha}(g_{\alpha, i})$
with
$a_{\alpha, i} \in R^*_{\alpha, h_{\alpha, i}}$
and the fact that
$p_{\alpha, i}$ does not divide $h_{\alpha, i}$
prove that the class
$\bar a_{\alpha}$
is unramified at each the
$p_{\alpha, i}$.
Thus  the class
$\bar a_{\alpha}$
is  indeed $R_{\alpha}$-unramified.

By purity for $R_{\alpha}$ there exists an
$a_{\alpha}^{\prime} \in R^*_{\alpha}$
such that
$\bar a_{\alpha}^{\prime}=\bar a_{\alpha}$ in
$\Cal F(K_{\alpha})$.
The exact sequence
$
1@>>>U_{A,\sigma}@>>>{\text{Sim}}_{A,\sigma}@>\mu>>{\Bbb G}_m@>>>1
$
of algebraic group schemes over $S_{\alpha}$ shows that the kernel
of the boundary map
$\partial: \Cal F(S_{\alpha}) \to H^1(S_{\alpha}, U_{A,\sigma})$
is trivial.
The Main Theorem of [Oj-P1] states that the kernel
$$
\ker[H^1(S_{\alpha}, U_{A,\sigma}) \to H^1(K_{\alpha}, U_{A,\sigma})]
$$
is trivial. Thus
$\Cal F(S_{\alpha})$ injects into $\Cal F(K_{\alpha})$ and
$\bar a_{\alpha}^{\prime}=\bar a_{\alpha}$
already in
$\Cal F(S_{\alpha})$.
The commutative diagram
$$
\CD
R_{\alpha} @>{\phi_{\alpha}}>> R \\
@VVV @VVV \\
S_{\alpha} @>>> K\ . \\
\endCD
$$
shows that
$\overline {\phi_{\alpha}(\alpha^{\prime})}=\bar a$
in $\Cal F(K)$.
This completes the proof of Theorem 1.2.

\enddemo

\enddocument